\documentclass[11pt]{amsproc}

\usepackage{latexsym}
\usepackage{amsmath}
\usepackage{amsfonts}
\usepackage{amssymb}
\usepackage{amscd}

\numberwithin{equation}{section}

\newtheorem{proposition}{Proposition}[section]
\newtheorem{lemma}[proposition]{Lemma}
\newtheorem{theorem}[proposition]{Theorem}

\def\itheorem#1#2{\newtheorem{#1}[proposition]{#2}}

\itheorem{Remark}{Remark}

\def\Hom{\mathop{\rm Hom}\nolimits}
\def\Ext{\mathop{\rm Ext}\nolimits}

\def\Q{{\mathbb{Q}}}
\def\Z{{\mathbb{Z}}}
\def\N{{\mathbb{N}}}

\def\tc{{\mathcal{T}\mathcal{C}}}

\begin{document}
\title[p-rank of {$\bf Ext(A,B)$} for countable abelian groups $\bf A$ and $\bf B$]{\bf On the p-rank of {$\bf Ext(A,B)$} for countable abelian groups $\bf A$ and $\bf B$}

\author{Lutz Str\"ungmann}

\address{Faculty of Mathematics,
University of Duisburg--Essen, 45117 Essen, Germany}

\email{lutz.struengmann@uni-due.de}
\author{Saharon Shelah}

\address{Department of Mathematics, The Hebrew University of Jerusalem, Israel
and Rutgers University, New Brunswick, NJ U.S.A.}

\email{shelah@math.huji.ac.il}

\thanks{This is number 968 in the first author's list of publication. The authors would like to thank the German Israeli Foundation for their support project No. I-963-98.6/2007 ; the second author was also supported by project STR 627/1-6 of the German Research Foundation DFG}

\maketitle

\begin{abstract}
In this note we show that the p-rank of $\Ext(A,B)$ for countable torsion-free abelian groups $A$ and $B$ is either countable or the size of the continuum.
\end{abstract}

\section{Introduction}
The structure of $\Ext(A,B)$ for torsion-free abelian groups $A$ has received much attention in the literature. In particular in the case of $B=\Z$ complete characterizations are available in various models of ZFC (see \cite[Sections on the structure of Ext]{EkMe}, \cite{ShSt1}, and \cite{ShSt2} for references). It is easy to see that $\Ext(A,B)$ is always a divisible group for
any torsion-free group $A$. Hence it is of the form
\[ \Ext(A,B)= \bigoplus\limits_{p \in
\Pi}\Z(p^{\infty})^{(\nu_p)} \oplus \Q^{(\nu_0)} \] for some
cardinals $\nu_p,\nu_0$ ($p \in \Pi)$ which are uniquely
determined. Here $\Z (p^{\infty})$ is the $p$-Pr\"ufer group and $\Q$ is the group of rational numbers. Thus, the obvious question that arises is which
sequences $(\nu_0, \nu_p : p \in \Pi)$ can appear as the cardinal
invariants of $\Ext(A,B)$ for some (which) torsion-free abelian
group $A$ and arbitrary abelian group $B$? As mentioned above for $B=\Z$ the answer is pretty much known but for general $B$
there is little known so far. Some results were obtained in \cite{Fr} and \cite{FrSt} for countable abelian groups $A$ and $B$. However, one essential question was left open, namely if the situation is similar to the case $B=\Z$ when it comes to $p$-ranks. It was conjectured that the $p$-rank of $\Ext(A,B)$ can only be either countable or the size of the continuum whenever $A$ and $B$ are countable. Here we prove that the conjecture is true.\\

Our notation is standard and we write maps from the left. If $H$ is a pure subgroup of the abelian group $G$, then we will write $H \subseteq_* G$. The set of natural primes is denoted by $\Pi$. For further details on abelian groups  and set-thoretic methods we refer to \cite{Fu} and \cite{EkMe}.

\section{Proof of the conjecture on the p-rank of Ext}

It is well-known that for torsion-free abelian groups $A$ the group of extensions $\Ext(A,B)$ is divisible for any abelian group $B$. Hence it is of the form
\[ \Ext(A,B)= \bigoplus\limits_{p \in
\Pi}\Z(p^{\infty})^{(\nu_p)} \oplus \Q^{(\nu_0)} \] for some
cardinals $\nu_p,\nu_0$ ($p \in \Pi)$ which are uniquely
determined.\\
 The invariant $r_p(\Ext(A,B)):=\nu_p$ is called the {\it $p$-rank} of $\Ext(A,B)$ while $r_0(\Ext(A,B)):=\nu_0$ is called the {\it torsion-free rank} of $\Ext(A,B)$. The following was shown in \cite{FrSt} and gives an almost complete description of the structure of $\Ext(A,B)$ for countable torsion-free $A$ and $B$. Recall that the {\it nucleus} $G_0$ of a torsion-free abelian group $G$ is the largest subring $R$ of $\Q$ such that $G$ is a module over $R$.

 \begin{proposition}
 Let $A$ and $B$ be torsion-free groups with $A$ countable and $\left|B\right| < 2^{\aleph_0}$. Then either
\begin{enumerate}
\item[{\bf i)}] $r_0(\Ext(A,B))=0$ and $A \otimes B_0$ is a free $B_0$-module \underline{or}
\item[{\bf ii)}] $r_0(\Ext(A,B))=2^{\aleph_0}$.
\end{enumerate}
\end{proposition}

Recall that for a torsion-free abelian group $G$ the {\it $p$-rank} $r_p(G)$ is defined to be the $\Z/p\Z$-dimension of the vectorspace $G/pG$.

\begin{proposition} The following holds true:
\begin{enumerate}
\item[{\bf i)}] If $A$ and $B$ are countable torsion-free abelian groups and $A$ has finite rank then $r_p(\Ext(A,B)) \leq \aleph_0$. Moreover, all cardinals $\leq \aleph_0$ can appear as the $p$-rank of some $\Ext(A,B)$.
\item[{\bf ii)}] If $A$ and $B$ are countable torsion-free abelian groups with $\Hom(A,B)=0$ then $$r_p(\Ext(A,B))=\left\{ \begin{array}{c c} 0 & if \ r_p(A)=0 \ or \ r_p(B)=0 \\ r_p(A) \cdot r_p(B) & if \ 0<r_p(A),r_p(B) < \aleph_0 \\ \aleph_0 & if \ 0<r_p(A)<\aleph_0 \ and \ r_p(B)=\aleph_0 \\ 2^{\aleph_0} & if \ r_p(A)=\aleph_0 \ and \ 0<r_p(B) \leq \aleph_0 \end{array} \right.$$.
\end{enumerate}
 \end{proposition}

Recall the following way of calculating the p-rank of $\Ext(A,B)$ for torsion-free groups $A$ and $B$ (see \cite[page 389]{EkMe}).

\begin{lemma} \label{rechnen}
For a torsion-free abelian group $A$ and an arbitrary abelian group $B$ let $\varphi^{p}$ be the map that sends $\psi \in \Hom(A,B)$ to $\pi \circ \psi$ with $\pi \ : \ \ B \ \longrightarrow \ B/pB$ the canonical epimorphism. Then
\[r_p(\Ext(A,B))=dim_{\Z/p\Z}(\Hom(A,B/pB)/\varphi^{p}(\Hom(A,B))).\]
\end{lemma}

We now need some preparation from descriptive set-theory in order to prove our main result. Recall that a subset of a topological space is called {\it perfect} if it is closed and contains no isolated points. Moreover, a subset $X$ of a Polish space $V$ is {\it analytic} if there is a Polish space $Y$ and a Borel (or closed) set $B \subseteq V \times Y$ such that $X$ is the projection of $B$; that is,
\[ X=\{v\in V|(\exists y\in Y)\langle v,y \rangle\in B\}. \]

\begin{proposition}\label{analytic}
If $E$ is an analytic equivalence relation on $\Gamma=\{ X : X \subseteq \omega\}$ that satisfies
\[ (\dag) \quad \textit{ if } X,Y \subseteq \omega, n \not\in Y, X=Y \cup \{n\} \textit{ then $X$ and $Y$ are not $E$-equivalent } \]
then there is a perfect subset $T$ of $\Gamma$ of pairwise nonequivalent $X \subseteq \omega$.
\end{proposition}

\proof See \cite[Lemma 13]{Sh2}. \qed\\

We are now in the position to prove our main result.

\begin{theorem}
Let $A$ be a countable torsion-free abelian group and $B$ an arbitrary countable abelian group. Then either
\begin{itemize}
\item $r_p(\Ext(A,B)) \leq \aleph_0$ \underline{or}
\item $r_p(\Ext(A,B))=2^{\aleph_0}$.
\end{itemize}
\end{theorem}

\proof The proof uses descriptive set theory, relies on \cite[Lemma 13]{Sh2} and is inspired by \cite[Lemma 2.2]{HaSh}. By the above Lemma \ref{rechnen} the p-rank of $\Ext(A,B)$ is the dimension $\kappa$ of the $\Z/p\Z$ vector space $L:=\Hom(A,B/pB)/\varphi^{p}(\Hom(A,B))$ where $\varphi^p$ is the natural map. We choose a basis $\left\{[\varphi_{\alpha}] \ | \ \alpha < \kappa \right\}$ of $L$ with $\varphi_{\alpha} \in \Hom(A,B/pB)$ and assume that $\aleph_0 < \kappa$. Note that $[\varphi_{\alpha}] \neq 0$ means, that $\varphi_{\alpha} \ : \ A \rightarrow B/pB$ has no lifting to an element $\psi \in \Hom(A,B)$ such that $\varphi_{\alpha}=\pi \circ \psi$.\\
Now let $A=\bigcup\limits_{i < \omega} A_i$ with $rk(A_i)$ finite. By a pigeonhole-principle there are $\alpha_1 \neq \beta_1$ such that $\varphi_{\alpha_1} \upharpoonright A_1=\varphi_{\beta_1} \upharpoonright A_1$. We define $\psi_1 := \varphi_{\alpha_1}-\varphi_{\beta_1}$ and obtain $A_1 \subseteq Ker(\psi_1)$. Obviously, $\psi_1$ has no lifting because $\left\{[\varphi_{\alpha}] \ | \ \alpha < \kappa \right\}$ is a basis of $L$, hence $[\psi_1] \not= 0$. Since $\alpha_1 \neq \beta_1$ there exists $x_1 \in A$ satisfying $\varphi_{\alpha_1}(x_1) \neq \varphi_{\beta_1}(x_1)$. Let $n$ be minimal with $x_1 \in A_n \setminus A_{n-1}$. Because $A_n$ has finite rank we similarly get $\alpha_2 \neq \beta_2$ such that $\varphi_{\alpha_2} \upharpoonright A_n=\varphi_{\beta_2} \upharpoonright A_n$ and define
$\psi_2 := \varphi_{\alpha_2}-\varphi_{\beta_2}$ with $A_n \subseteq Ker(\psi_2)$. Clearly we have $\psi_1 \neq \psi_2$ since $\psi_1(x_1) \neq 0 = \psi_2(x_1)$ and $[\psi_2]\not=0$. Continuing this construction we get $\aleph_0$ pairwise different $\psi_n \in \Hom(A,B/pB)$ with $[\psi_n]\not= 0$. \\
Now let $\eta \in \ ^{\omega}2$, which means that $\eta$ is a countable $\left\{0,1\right\}$-sequence and choose $\psi_{\eta} := \sum\limits_{n \in supp(\eta)} \psi_n$ where $supp(\eta):=\left\{n \in \omega \ | \ \eta(n)=1\right\}$. This is well-defined since for any $a \in A$, the sum $\psi_{\eta}(a)$ consists of only finitely many summands because there is $n \in \N$ such that $a \in A_n$ and $A_n$ is in the kernel of $\psi_m$ for sufficiently large $m$. Note that also $\psi_{\eta} \neq \psi_{\eta^{\prime}}$ for all $\eta \neq \eta^{\prime}$. It remains to prove that the size of $\{ [\psi_{\eta}] : \eta \in \ ^{\omega}2\}$ is $2^{\aleph_0}$ which then implies that $\kappa=2^{\aleph_0}$.\\
We now define an equivalence relation on $\Gamma=\{X \subseteq \omega \}$ in the following way: $X \sim X'$ if and only if $[\psi_{\eta_X}]=[\psi_{\eta_{X'}}]$, where $\eta_X \in {}^{\omega}2$ is the characteristic function of $X$. We claim that\\

\centerline{($\ddag$) \quad for $X,X' \subseteq \omega$ and $n \not\in X$ such that $X'=X \cup \{n\}$ we have $X \not\sim X'$.}
\quad\\
But this is obvious since in this case $\psi_{\eta_X} - \psi_{\eta_{X'}}=-\psi_n$ and $[-\psi_n]\not=0$, so $[\psi_{\eta_X}] \not= [\psi_{\eta_{X'}}]$. We claim that the above equivalence relation is analytic. It then follows from Proposition \ref{analytic} that there is a perfect subset of $\Gamma$ of pairwise nonequivalent $X \subseteq \omega$ and hence there are $2^{\aleph_0}$ distinct equivalence classes for $\sim$. Thus the $p$-rank of $\Ext(A,B)$ must be the size of the continuum as claimed.\\
In order to see that $\sim$ is analytic recall that both, the cantor space ${}^{\omega}2$ (and hence $\Gamma$) and the Baire space ${}^{\omega}\omega$ are Polish spaces with the natural topologies. For $\sim$ to be analytic we therefore have to show that $\Delta=\{(X,X') : X,X' \in \Gamma \textit{ and } X \sim X' \}$ is an analytic subset of the product space $\Gamma \times \Gamma$. We enumerate $A=\{ a_n : n \in \omega\}$ and $B=\{ b_n : n \in \omega\}$ and let $C_1 \subseteq {}^{\omega}\omega$ be the set of all $\rho \in {}^{\omega}\omega$ such that $\rho$ induces a homomorphism $f_{\rho} \in \Hom(A,B/pB)$, i.e. the map $f_{\rho} : A \rightarrow B/pB$ sending $a_n$ onto $b_{\rho(n)} + pB$ is a homomorphism. $C_2$ is defined similarly replacing $\Hom(A,B/pB)$ by $\Hom(A,B)$, i.e. $\nu \in C_2$ if the map $g_{\nu}: A \rightarrow B$ with $a_n \mapsto b_{\nu(n)}$ is a homomorphism. Clearly $C_1$ and $C_2$ are closed subsets of ${}^{\omega}\omega$. Put $P=\{ (\rho_1,\rho_2,\nu) : \rho_1,\rho_2 \in C_1; \nu \in C_2 \text{ and } f_{\rho_1}-f_{\rho_2} = \varphi^p (g_{\nu})\}$. Then $P$ is a closed subset of ${}^{\omega}\omega \times {}^{\omega}\omega \times {}^{\omega}\omega$. Now the construction above gives a continuous function $\Theta$ from ${}^{\omega}2$ to ${}^{\omega}\omega$ by sending $\eta$ to $\rho=\Theta(\eta)$ where $\rho$ is induced by the homomorphism $\psi_{\eta}$. Then $X \sim X'$ if and only if there is $\nu \in C_2$ such that $(\Theta(\eta_X),\Theta(\eta_{X'}),\nu) \in P$ and thus $\sim$ is an analytic equivalence relation.
\qed


{\small
\begin{thebibliography}{99}
\addcontentsline{toc}{chapter}{Literaturverzeichnis}


\bibitem[EkMe]{EkMe}{\bf P.C. Eklof and A. Mekler}, {\it Almost Free Modules, Set-
Theoretic Methods (revised edition)}, Amsterdam, New York,
North-Holland, Math. Library (2002).


\bibitem[Fr]{Fr}{\bf S. Friedenberg}, {\it Extensions of Abelian Groups}, S\"udwestdeutscher Verlag f\"ur Hochschulschrift (2010).

\bibitem[FrSt]{FrSt}{\bf S. Friedenberg and L. Str\"ungmann}, {\it The structure of $\Ext(A,B)$ for countable abelian groups $A$ and $B$}, submitted (2010).

\bibitem[Fu]{Fu}{\bf L. Fuchs}, {\it Infinite Abelian Groups, Vol. I and II},
Academic Press (1970 and 1973).



\bibitem[HaSh]{HaSh}{\bf L. Harrington and S. Shelah}, {\it Counting equivalence classes for co-$\kappa$-Souslin relation}, Proc. Conf. Prague, 1980, Logic Colloq. North-Holland, 1982, 147-152.











\bibitem[Sh]{Sh2}{\bf S. Shelah}, {\it Can the fundamental (homotopy) group of a space be the rationals?}, Proc. AMS {\bf 103} (1988),
627--632.

\bibitem[ShSt1]{ShSt1}{\bf S. Shelah and L. Str\"ungmann}, {\it A characterization of $\Ext(G,\Z)$ assuming (V=L)}, Fundamenta Math. {\bf 193} (2007), 141-151.

\bibitem[ShSt2]{ShSt2}{\bf S. Shelah and L. Str\"ungmann}, {\it On the p-rank of $\Ext(G,\Z)$ in certain models of ZFC}, Algebra and Logic {\bf 46} (2007), 200-215.


\end{thebibliography}
}
\end{document}